\documentclass[reqno]{article}
\usepackage{comment,subeqnarray,amsthm,amsmath,amssymb}
\usepackage{hyperref}
\usepackage{epsfig,color,datetime}
\usepackage{multirow}
\usepackage{colortbl,array}

\usepackage[table]{xcolor}
\usepackage[numbers, square, comma]{natbib}
\usepackage{algpseudocode,algorithm,algorithmicx}
\usepackage{setspace}
\usepackage{graphicx}

\numberwithin{equation}{section}
\numberwithin{figure}{section}
\numberwithin{table}{section}
\newtheorem{theorem}{Theorem}[section]

\newtheorem{proposition}[theorem]{Proposition}

\newtheorem{example}[theorem]{Example}

\newcommand{\exref}[1]{Example~\ref{#1}}

\newcommand{\eq}[1]{\begin{eqnarray}\label{#1}}
\newcommand{\qe}{\end{eqnarray}}

\newcommand{\be}{\begin{eqnarray}}
\newcommand{\ee}{\end{eqnarray}}
\newcommand{\bal}{\begin{aligned}}
\newcommand{\eal}{\end{aligned}}
\newcommand{\bes}{\begin{eqnarray*}}
\newcommand{\ees}{\end{eqnarray*}}
\newcommand{\bs}{\begin{subeqnarray}}
\newcommand{\es}{\end{subeqnarray}}
\newcommand{\bss}{\begin{subeqnarray*}}
\newcommand{\ess}{\end{subeqnarray*}}

\newcounter{saveeqn}

\def\rand{\operatorname{rand}}

\def\p{\partial}

\def\O{\Omega}

\def\NChz{{\mathcal NC}^h_0}
\def\NChz2d{{[\mathcal{NC}}^h_0]^2}
\def\tNChz2d{\widetilde {[\mathcal{NC}}^h_0]^2}

\def\and{\quad\text{and}\quad}

\def\<{\left\langle}
\def\>{\right\rangle}

\def\mbA{\mathbf A}
\def\mbb{\mathbf b}
\def\mbeta{\boldsymbol \eta}

\def\b1{\mathbf 1}

\def\Tau{{\mathcal T}}

\newcolumntype{x}[1]{>{\centering\hspace{0pt}}p{#1}} 

\newcommand{\vertiii}[1]{{\left\vert\kern-0.25ex\left\vert\kern-0.25ex\left\vert #1 \right\vert\kern-0.25ex\right\vert\kern-0.25ex\right\vert}}

\def\balpha{\boldsymbol\alpha}
\def\bgamma{\boldsymbol\gamma}
\def\boldeta{\boldsymbol\eta}
\def\O{\Omega}

\def\Tau{\mathcal T}

\begin{document}

\allowdisplaybreaks

\newif\iflong
\longfalse

\title{Algebraic Multiscale Method for one--dimensional elliptic problems}

\author{
Kanghun Cho%
\thanks{Samsung Fire \& Marine Insurance Co., Ltd., 14, Seocho-daero 74-gil, Seocho-gu, Seoul 06620,
         Korea;}
,$\quad$ Roktaek Lim%
\thanks{Department of Biology, Hong Kong Baptist University,
  Cha Chi-ming Science Tower, Ho Sin Hang Campus,
Kowloon Tong, Hong Kong}
,$\quad$ Dongwoo Sheen%
\thanks{Department of Mathematics, Seoul National University,
  Seoul 08826, Korea
  Eemails:serein@snu.ac.kr, rokt.lim@gmail.com, sheen@snu.ac.kr}}

\maketitle

\begin{abstract}
      In this paper we propose an idea of constructing a macro--scale matrix
system
given a micro--scale matrix linear system. Then the
macro--scale system is solved at cheaper computing costs.
The method uses the idea of the generalized multiscale finite element method based. Some numerical results are presented.
\end{abstract}


%
%

\section{Introduction}
In this paper we propose an algebraic multiscale method for
one--dimensional elliptic problems. The extension to two dimensional
case will appear in \cite{cho2022algebraic2d}.

Assume that only the algebraic information on the components of
a micro--scale linear system are known, but no further information on the
coefficient $\kappa$ and the source term $f$ are available.
In this situation, our object is to try to construct macro--scale
linear systems using accessible information and find numerical
solutions which possess similar properties of the solutions
obtained by multiscale methods. In some sense, this is an inverse
problem to fetch the necessary information on the
fast--varying coefficient and the source function of governing
elliptic equation. This will give us elliptic equation in micro--scale.
Then, we follow the standard approach to build a macro--scale linear
system based on a multiscale method. Although other multiscale methods
can generate similar macro--scale linear systems, in this paper we use
the generalized multi-scale finite element method approach.

 Multiscale methods have been actively developed in various manners including heterogeneous multiscale methods \cite{abdulle2009finite, abdulle2012heterogeneous}, multiscale hybridizable discontinuous Galerkin methods \cite{cho2020multiscale, efendiev2015spectral, efendiev2015multiscale}, and multiscale finite element methods \cite{efendiev2009multiscale,
efendiev2000convergence,hou1999convergence}. Here we adopt the generalized multiscale finite element method(GMsFEM) \cite{efendiev2013generalized, lee2017nonconforming}.

The generalized multiscale finite element spaces consist of snapshot function spaces, offline function spaces, and moment function spaces. First, snapshot functions are obtained by solving $\kappa$--harmonic problems in each macro element. Then offline functions are constructed by applying suitable dimension reduction techniques to snapshot function space. We choose offline functions identical to snapshot function space since there are only two snapshot functions in each macro element. 

This paper is organized as follows. In section 2, we briefly review the nonconforming generalized multiscale finite element method(GMsFEM) based on DSSY finite element space. Then the algebraic multiscale method for two--dimensional elliptic problem is introduced in section 3 following the framework of GMsFEM. Section 4 is devoted to energy norm error estimate of the proposed method.
In Section 5, representative numerical results are presented. A conclusion is given in Section 6.

%
%
\section{Preliminaries}
In this paper, we consider the following one-dimensional elliptic problem:
\begin{equation}\label{eq:1d-ell}
\left\{
\begin{aligned}
-\frac{d}{dx} \big( A \frac{d u}{dx} \big) = f & \text{ in } \O,\\
u = 0, & \text{ if }  \text{ on } \p\O,
\end{aligned}
\right.
\end{equation}
where $\O=(0,1),\p\O=\{0,1\},$ and $A$ is a rapidly varying coefficient. Denote by $(\Tau_H)_{0<H<1}$ and $(\Tau_h)_{0<h<1}$ two families of
macro and micro--scale triangulations of
$(0,1)$ into macro and micro--scale subintervals such that $0=X^{0}< X^1 < \cdots,< X^{N_H}=1$ and
$0=x_{0}< x_1 < \cdots < x_{N_h}=1.$ Here, and in what follows,
$H$ and $h$ stand for the macro and micro--scale mesh parameters given by
\[
H = \max_{K=1,\cdots, N_H} (X^K - X^{K-1}),\quad
  h = \max_{j=1,\cdots, N_h} (x_j - x_{j-1}).
    \]
We assume that $\{X^0,X^1,\cdots,X^{N_H}\} \subset
\{x_0,x_1,\cdots,x_{N_h}\}$ and $h\ll H<1.$
For $K=1,\cdots, N_H,$ denote by $H^K$ the size of $K$-th macro interval $I^{K}=(X^{K-1},X^{K}).$
Let $\{ x^{K}_{j} \}_{j=0}^{N^K}$ be the set of nodes for the macro
interval $I^{K}=(X^{K-1},X^{K})$ and designate by $I^K_j$ the $j$--th subinterval
$(x^{K}_{j-1},x^{K}_{j})$ with length $h^K_j$ for $j = 1,\cdots, N^K$ such that
$x^{K}_{0}  = X^{K-1}$ and $x^{K}_{N^K}  = X^{K}.$
For each $K,$ let $\{ \phi^{K}_{j} \}_{j=0,\ldots,N^K}$ be the space of standard basis functions for the $C^0$--piecewise linear finite element space on the interval $I^{K}=(X^{K-1},X^{K}).$
Denote by $\Psi^{K}_{\pm}$ the macro--scale basis function
in the interval $I^{K},$ which can be obtained as the solutions of 
\begin{equation}\label{eq:macrobasis-}
\left\{
\begin{aligned}
&-\frac{d}{dx} \big( A \frac{d}{dx} \Psi^{K}_{-} \big) = 0 \text{ in } I^{K},\\
&\Psi^{K}_{-}(X^{K-1})=1,\; \Psi^{K}_{-}(X^{K})=0,
\end{aligned}
\right.
\end{equation}
and
\begin{equation}\label{eq:macrobasis+}
\left\{
\begin{aligned}
&-\frac{d}{dx} \big( A \frac{d}{dx} \Psi^{K}_{+} \big) = 0 \text{ in } I^{K},\\
&\Psi^{K}_{+}(X^{K-1})=0,\; \Psi^{K}_{+}(X^{K})=1.
\end{aligned}
\right.
\end{equation}

For each $K$ let us seek $\Psi^{K}_{h,\pm}$ which approximate $\Psi^{K}_{\pm}$
in the form
\begin{subeqnarray}\label{eq:PsiKhpm}
\Psi^{K}_{h,-} &=& \sum_{j=1}^{N^K-1} \eta^K_{j,-} \phi^{K}_{j} + \phi^{K}_{0},\\
\Psi^{K}_{h,+} &=& \sum_{j=1}^{N^K-1} \eta^K_{j,+} \phi^{K}_{j} + \phi^{K}_{N^K}.
\end{subeqnarray}

Since $\Psi^{K}_{h,\pm}$ are piecewise-linear in $I^{K}$, one may set
\begin{equation}\label{eq:diffPsiKhpm}
\frac{d}{dx} \Psi^{K}_{h,\pm} = \gamma^{K}_{j,\pm} \; \text{in} \;
I^{K}_j\quad\text{for some constant }\gamma^{K}_{j,\pm},j=1,\cdots, N^K.
\end{equation}

Assuming that $\frac1{A(x)} \in L^1{(a,b)},$ one sees that the exact solution of
the differential equation
\begin{equation*}
\left\{
\begin{aligned}
&-\frac{d}{dx} \big( A \frac{dw}{dx} \big) = 0, \; \text{ in } (a,b)\\
&w(a) = 0,\; w(b)=1.
\end{aligned}
\right.
\end{equation*}
is given by
\begin{equation}\label{eq:exact-sol-w}
w(x) = \beta \int_{a}^{x} \frac{1}{A(s)} \; ds \quad a.e.\,
 x \in (a,b),
\end{equation}
where $\beta$ is the harmonic mean of $A(x)$ over $(a,b),$ {\it i.e.,}
\begin{equation}\label{eq:exact-sol-beta}
\beta = \frac{1}{ \int_{a}^{b} \frac{1}{A(s)} \; ds } \;.
\end{equation}
%
Thanks to \eqref{eq:exact-sol-w} and \eqref{eq:exact-sol-beta},
it is easy to see that
\begin{subeqnarray}\label{eq:formularPsiKhpm}
  \Psi^K_-(x) &=&
 \frac{\int^{X^{K}}_x \frac{1}{A(s)} \; ds }{\int_{X^{K-1}}^{X^K}
  \frac{1}{A(s)} \; ds}
  \text{ for all } x \in I^K,\\
  \Psi^K_+(x) &=&
 \frac{\int_{X^{K-1}}^x \frac{1}{A(s)} \; ds }{\int_{X^{K-1}}^{X^K}
  \frac{1}{A(s)} \; ds}
  \text{ for all } x \in I^K.
\end{subeqnarray}
Denote an $n$--dimensional vector with parameters $K$ and $\pm$ as follows:
\[
\balpha^{K}_{\pm} = (\alpha^{K}_{1,\pm},\cdots, \alpha^{K}_{n,\pm})^t \in \mathbb R^{n}.
\]
Recalling \eqref{eq:PsiKhpm}, \eqref{eq:diffPsiKhpm}, and
\eqref{eq:formularPsiKhpm}, and utilizing the principle of energy norm minimization of finite
element method,
we deduce the following equalities:
\begin{equation*}
\begin{aligned}
&\min_{\boldeta^K_{+}\in \mathbb R^{N^K}}
\Big\{
\int_{I^{K}} A(x) \big[ \frac{d}{dx} \big( \Psi^{K}_+(x) - \Psi^{K}_{h,+}(x) \big) \big]^{2} \; dx
\Big\}^{\frac12} \\
&=
\min_{\boldeta^K_{+}\in \mathbb R^{N^K}}
\sum_{j=1}^{N^K}
\Big\{
\int_{I^{K}_{j}}
A(x) \big[ \frac{d}{dx} \big( \Psi^{K}_+(x) - \Psi^{K}_{h,+}(x) \big) \big]^{2} \; dx
\Big\}^{\frac12} \\
&=
\min_{\bgamma^{K}_{+}} \sum_{j=1}^{N^K}
\Big\{
\int_{I^{K}_{j}}
A(x) \big( \frac{\beta}{A(x)} - \gamma^{K}_{j,+} \big)^{2} \; dx
\Big\}^{\frac12}\\
&=
\min_{\gamma^{K}_{j}} \sum_{j=1}^{N}
\Big\{
\int_{I^{K}_{j}}
\big( \frac{\beta^{2}}{A(x)} - 2\beta \gamma^{K}_{j,+} + A(x) (\gamma^{K}_{j,+})^{2} \big) \; dx
\Big\}^{\frac12}.
\end{aligned}
\end{equation*}
After differentiating the above with respect to $\gamma^{K}_{j}$, we have
\begin{equation*}
\begin{aligned}
\gamma^{K}_{j,+} = \frac{\int_{x^{K}_{j-1}}^{x^{K}_{j}} \beta \; dx}{\int_{x^{K}_{j-1}}^{x^{K}_{j}} A(x) \; dx}
= \frac{\beta h^{K}_{j}}{\int_{x^{K}_{j-1}}^{x^{K}_{j}} A(x) \; dx}.
\end{aligned} 
\end{equation*}
If $A(x)=1$, $\gamma^{K}_{j}=1$.

%
%

\section{Algebraic Multiscale Method}
In this section, we introduce an algebraic multiscale method.
Assume that we are given a linear system:
\begin{eqnarray}\label{eq:micro-lin-sys}
  \mbA^h \mbeta^h = \mbb^h,
\end{eqnarray}
which is obtained by a discretization of a micro--scale elliptic
equation of form \eqref{eq:1d-ell}. Assuming that
the coefficient $A$ and the source function $f$ are not known directly,
our aim is to construct a macro--scale matrix system
\begin{eqnarray}\label{eq:macro-lin-sys}
  \mbA^H \mbeta^H = \mbb^H
\end{eqnarray}
from the micro--scale linear system \eqref{eq:micro-lin-sys}.
In particular, without solving the macro--scale basis problem
\eqref{eq:macrobasis-} for each macro element $I^K,$ we try to infer
the components of $\mbA^H$ and $\mbb^H$ which are obtained
by a multiscale method from the
structure of the elliptic problem \eqref{eq:1d-ell}. Then the
macro--scale linear system \eqref{eq:macro-lin-sys} is solved at a
cheaper cost.

From now on, the midpoint rule is assumed to approximate integrals.

We state the algorithm as follows, and describe the details in the
subsections to follow.
\begin{enumerate}
\item[Step 1.] Approximate the coefficients and RHS of \eqref{eq:1d-ell} from the micro--scale matrix system;
\item[Step 2.] Construct a macro--scale matrix system from the
  information obtained in Step 1;
\item[Step 3.] Solve the macro--scale matrix system to get a multiscale solution.
\end{enumerate}
\subsection{Micro--scale problem}
Let $\{ \phi_{j} \}_{j=0,\ldots,N_h}$ be the space of standard basis functions for the $C^0$--piecewise linear finite element space on $\Omega=(0,1).$ We have
\begin{equation}\label{eq:micro_grad}
\frac{d \phi_{j}}{dx} \Big|_{ [x_{j-1},x_{j}] } = \frac{1}{h_j} \; \text{and} \;
\frac{d \phi_{j}}{dx} \Big|_{ [x_{j},x_{j+1}] } = -\frac{1}{h_{j+1}},
\end{equation}
where $h_{j}$ is the size of $j-$th micro interval $I_{j}=(x_{j-1},x_{j}).$
Let $A^h$ be the micro--scale stiffness matrix. For $1 \leq j \leq N_{h}-1,$ the diagonal element of $A^h$ is given as follows:
\begin{eqnarray}
\big[ A^h \big]_{j,j} &=& \int_{x_{j-1}}^{x_{j+1}} A(x) \frac{d
                          \phi_{j}}{dx} \frac{d \phi_{j}}{dx} \;
                          dx\nonumber \\
&=& \int_{x_{j-1}}^{x_{j}} A(x) (\frac{1}{h_{j}})^{2} \; dx
+ \int_{x_{j}}^{x_{j+1}} A(x) (-\frac{1}{h_{j+1}})^{2} \; dx\nonumber \\
&\approx& \frac{1}{h_{j}} A_{j-\frac12} +
          \frac{1}{h_{j+1}}A_{j+\frac12}.
          \label{eq:micro_stiff_diag}
\end{eqnarray}
For $j=0$ and $j=N_h,$
\begin{equation}\label{eq:micro_stiff_diag2}
\quad \big[ A^h \big]_{0,0} \approx \frac{1}{h_{1}}A_{\frac12}, \quad \big[ A^h \big]_{N_{h},N_{h}} \approx \frac{1}{h_{N_{h}}}A_{N_{h}-\frac12}.
\end{equation}
The off--diagonal element of $A^h$ can be computed similarly:
\begin{equation}\label{eq:micro_stiff_off-diag}
\big[ A^h \big]_{j,j-1} \approx -\frac{1}{h_{j}}A_{j-\frac12}, \quad \big[ A^h \big]_{j,j+1} \approx= -\frac{1}{h_{j+1}}A_{j+\frac12}.
\end{equation}
From above equation, we know that there is an one-to-one correspondence between the off--diagonal element of $A^h$ and the average of coefficient matrix in each micro interval.

To compute the right hand side $b^{h}$, we assume $f=\frac{dg}{dx}$ for some $g \in H^1(\Omega).$ Since the derivatives of micro--scale basis functions are piecewise constant functions, it makes the computation simpler. For $1 \leq j \leq N_{h}-1,$
\begin{equation}\label{eq:rhs_micro}
\begin{aligned}
b^{h}_{j} &= \int_{x_{j-1}}^{x_{j+1}} f \phi_{j} \; dx = \int_{x_{j-1}}^{x_{j}} \frac{d g}{dx} \phi_{j} \; dx + \int_{x_{j}}^{x_{j+1}} \frac{d g}{dx} \phi_{j} \; dx  \\
&= \big[ \phi_{j} g \big]_{x_{j-1}}^{x_{j}} - \int_{x_{j-1}}^{x_{j}} \frac{d \phi_j}{dx} g \; dx   + \big[ \phi_{j} g \big]_{x_{j}}^{x_{j+1}} - \int_{x_{j}}^{x_{j+1}} \frac{d \phi_j}{dx} g \; dx \\
&= g(x_j) -  \int_{x_{j-1}}^{x_{j}} \frac{1}{h_{j}} g \; dx - g(x_j) +  \int_{x_{j}}^{x_{j+1}} \frac{1}{h_{j+1}} g \; dx \\
&\approx g(x_{j+\frac12}) - g(x_{j-\frac12}).
\end{aligned}
\end{equation}
For $j=0$ and $j=N_h$,
\begin{equation}\label{eq:rhs_micro2}
b^{h}_{0} = g(x_{\frac12}) - g(x_0), \quad b^{h}_{N_{h}}=g(x_{N_{h}})-g(x_{N_{h}-\frac12}).
\end{equation}
Under the additional assumption $f=\frac{dg}{dx}$, there are only
$N_{h}+1$ equations \eqref{eq:rhs_micro}, \eqref{eq:rhs_micro2} to
decide $N_{h}+2$ unknowns $g(x_{0}),g(x_{\frac12}),g(x_{\frac32}),
\cdots, g(x_{N_{h}-\frac12}), g(x_{N_{h}}).$ Thus these unknowns are
computed up to an additive constant, which corresponds to the constant
of indefinite integration. Hence we further assume that
\begin{equation}\label{eq:rhs_micro3}
g(x_{0})=0,
\end{equation}
From \eqref{eq:rhs_micro}--\eqref{eq:rhs_micro3} it follows that
\begin{eqnarray}\label{eq:rhs_micro4}
  g(x_{0})=0;\, 
  g(x_{j+\frac12})=\sum_{k=0}^j b_k^h, j=0,\cdots,N_h-1,\;
  g(x_{N_h})=\sum_{k=0}^{N_h} b_k^h.
\end{eqnarray}

\subsection{Macro--scale problem}
Let $\{ \Psi^{K} \}_{I=0}^{N_{H}}$ be a set of macro--scale basis functions on $\Omega=(0,1).$ For $1 \leq K \leq N_{H}-1,$ $\Psi^{K}$ satisfies the equations:
\begin{equation}\label{eq:macro-basis-problem1}
\left\{
\begin{aligned}
&-\frac{d}{dx} \big( A \frac{d}{dx} \Psi^{K} \big) = 0 \text{ in } I^{K},\\
&\Psi^{K}(X^{K-1})=0,\; \Psi^{K}(X^{K})=1,
\end{aligned}
\right.
\end{equation}
and
\begin{equation}\label{eq:macro-basis-problem2}
\left\{
\begin{aligned}
&-\frac{d}{dx} \big( A \frac{d}{dx} \Psi^{K} \big) = 0 \text{ in } I^{K+1},\\
&\Psi^{K}(X^{K})=1,\; \Psi^{K}(X^{K+1})=0.
\end{aligned}
\right.
\end{equation}
Thus
\begin{equation*}
A \frac{d}{dx} \Psi^{K}= \left\{
\begin{aligned} &c^{K}  \quad \quad a.e. \text{ in } I^{K}, \\
				&c^{K+1}\quad a.e. \text{ in } I^{K+1},
\end{aligned}
\right.
\end{equation*}
where $c^{K}$ and $c^{K+1}$ are constants. Assume that $\frac{1}{A(x)} \in L^1{(a,b)}.$ Integrating $\frac{d}{dx} \Psi^{K} = c^{K} /A(x)$ over $I^{K}$ gives
\begin{equation*}
1 = c^{K}{\int_{X^{K-1}}^{X^{K}} \frac{dx}{A(x)}}.
\end{equation*}
That is,
\begin{equation*}
A \frac{d}{dx} \Psi^{K} = \frac{1}{{\int_{X^{K-1}}^{X^{K}} \frac{dx}{A(x)}}} \text{ on } I^{K}.
\end{equation*}
Similarly,
\begin{equation*}
A \frac{d}{dx} \Psi^{K} = -\frac{1}{{\int_{X^{K}}^{X^{K+1}} \frac{dx}{A(x)}}} \text{ on } I^{K+1}.
\end{equation*}
Let $A^H$ be the macro--scale stiffness matrix. For $1 \leq K \leq N_{H}-1,$ the diagonal element of $A^H$ is given as follows:
\begin{eqnarray}
\big[ A^H \big]_{K,K} &=& \int_{X^{K-1}}^{X^{K+1}} A(x) \frac{d \Psi^{K}}{dx} \frac{d \Psi^{K}}{dx} \; dx\nonumber\\
&=&\frac{1}{{\int_{X^{K-1}}^{X^{K}} \frac{dx}{A(x)}}} \int_{X^{K-1}}^{X^{K}} \frac{d \Psi^{K}}{dx} \; dx -  \frac{1}{{\int_{X^{K}}^{X^{K+1}} \frac{dx}{A(x)}}} \int_{X^{K}}^{X^{K+1}} \frac{d \Psi^{K}}{dx} \; dx\nonumber \\
&=& \frac{1}{{\int_{X^{K-1}}^{X^{K}} \frac{dx}{A(x)}}} +
    \frac{1}{{\int_{X^{K}}^{X^{K+1}} \frac{dx}{A(x)}}}.
    \label{eq:macro_stiff_diag}
\end{eqnarray}
For $K=0$ and $K=N_{H},$
\begin{equation}\label{eq:macro_stiff_diag2}
\big[ A^H \big]_{0,0} = \frac{1}{{\int_{X^{0}}^{X^{1}} \frac{dx}{A(x)}}}, \quad \big[ A^H \big]_{N_{H},N_{H}} = \frac{1}{{\int_{X^{N_{H-1}}}^{X^{N_{H}}} \frac{dx}{A(x)}}}.
\end{equation}

The off--diagonal element of $A^H$ can be computed similarly:
\begin{equation}\label{eq:macro_stiff_off-diag}
\big[ A^H \big]_{K,K-1} = -\frac{1}{{\int_{X^{K-1}}^{X^{K}} \frac{dx}{A(x)}}}, \quad \big[ A^H \big]_{K,K+1} = -\frac{1}{{\int_{X^{K}}^{X^{K+1}} \frac{dx}{A(x)}}}.
\end{equation}

\subsection{Algebraic formulation of macro--scale system}
We can construct a macro--scale matrix system from the micro--scale
matrix system.
Denote by $M_{K}$ the number of fine nodes on $(0,X^K).$ Then
\begin{eqnarray}\label{eq:integral_Ainv}
\int_{X^{K-1}}^{X^{K}} \frac{dx}{A(x)}  &=& \sum_{j=1}^{N^K} \int_{x^{K}_{j-1}}^{x^{K}_{j}} \frac{dx}{A(x)}
                                            \approx \sum_{j=1}^{N^K} \frac{h_{j}}{A(x^{K}_{j-\frac12})}
\\
                                            &=& \sum_{j=1}^{N^K} \frac{h_{j}}{A_{M_{K-1}+j-\frac12}}
=-\sum_{j=1}^{N^K} \Big(\big[A^h \big]_{M_{K-1}+j,M_{K-1}+j-1}
                                                \Big)^{-1}.
                                                \nonumber
\end{eqnarray}
Thus the off--diagonal element of $A^h$ can be computed by
\begin{equation}\label{eq:macro_stiff_alg_off-diag}
\big[ A^H \big]_{K,K-1} = -\frac{1}{{\int_{X^{K-1}}^{X^{K}} \frac{dx}{A(x)}}}
\approx \frac{1}{\sum_{j=1}^{N^K} \Big( \big[ A^h \big]_{M_{K-1}+j,M_{K-1}+j-1}\Big)^{-1} },
\end{equation}
and
\begin{equation}\label{eq:macro_stiff_alg_off-diag2}
\big[ A^H \big]_{K,K+1} = -\frac{1}{{\int_{X^{K}}^{X^{K+1}} \frac{dx}{A(x)}}}
\approx \frac{1}{\sum_{j=1}^{N^{K+1}} \Big( \big[ A^h \big]_{M_{K}+j,M_{K}+j-1}\Big)^{-1} }.
\end{equation}
We can compute the diagonal element of $A^h$ by the same way.
For $1 \leq K \leq N_{H}-1,$
\begin{eqnarray}\label{eq:macro_stiff_alg_diag}
\big[ A^H \big]_{K,K} &=& \frac{1}{{\int_{X^{K-1}}^{X^{K}} \frac{dx}{A(x)}}} + \frac{1}{{\int_{X^{K}}^{X^{K+1}} \frac{dx}{A(x)}}} \\ &\approx& -\frac{1}{\sum_{j=1}^{N^K} \Big( \big[ A^h \big]_{M_{K-1}+j,M_{K-1}+j-1}\Big)^{-1} } -  \frac{1}{\sum_{j=1}^{N^{K+1}} \Big( \big[ A^h \big]_{M_{K}+j,M_{K}+j-1}\Big)^{-1} }.\nonumber
\end{eqnarray}
For $K=0$ and $K=N_H,$
\begin{eqnarray*}\label{eq:macro_stiff_alg_diag2}
\big[ A^H \big]_{0,0} &=& \frac{1}{{\int_{X^{0}}^{X^{1}} \frac{dx}{A(x)}}} \approx -\frac{1}{\sum_{j=1}^{N^0} \Big( \big[ A^h \big]_{j,j-1}\Big)^{-1} }, \\
\big[ A^H \big]_{N_{H},N_{H}} &=& \frac{1}{{\int_{X^{N_{H-1}}}^{X^{N_{H}}} \frac{dx}{A(x)}}} \approx  -\frac{1}{\sum_{j=1}^{N^{N^H}} \Big( \big[ A^h \big]_{M_{N_{H}-1}+j,M_{N_H-1}+j-1}\Big)^{-1} }.\nonumber
\end{eqnarray*}
Thus we can construct the macro--scale stiffness matrix only using the micro--scale stiffness matrix elements.

For the right hand side, note that
\begin{eqnarray*}
\int_{X^{K-1}}^{X^{K}} f \Psi^{K} \; dx &=& \sum_{j=1}^{N^K} \int_{x^{K}_{j-1}}^{x^{K}_{j}} \frac{d g}{dx} \Psi^{K} \; dx 
= \sum_{j=1}^{N^K} \Big( \big[\Psi^{K} g \big]_{x^{K}_{j-1}}^{x^{K}_{j}} - \int_{x^{K}_{j-1}}^{x^{K}_{j}} \frac{d \Psi^{K}}{dx} g \; dx \Big) \\
&=& g(x^{K}_{N}) - \sum_{j=1}^{N^K} \int_{x^{K}_{j-1}}^{x^{K}_{j}} \frac{1}{A(x)}\frac{1}{{\int_{X^{K-1}}^{X^{K}} \frac{dx}{A(x)}}} g \; dx \\
&\approx& g(x^{K}_{N}) - \frac{1}{{\int_{X^{K-1}}^{X^{K}} \frac{dx}{A(x)}}} \sum_{j=1}^{N^K} h_{j} \frac{g(x^{K}_{j-\frac12})}{A(x^{K}_{j-\frac12})}.
\end{eqnarray*}
Thus for $1 \leq K \leq N_{H}-1,$
\begin{equation*}
\begin{aligned}
b^{H}_{K} &= \int_{X^{K-1}}^{X^{K+1}} f \Psi^{K} \; dx = \int_{X^{K-1}}^{X^{K}} f \Psi^{K} \; dx + \int_{X^{K}}^{X^{K+1}} f \Psi^{K} \; dx \\
&\approx \Big(g(x^{K}_N) - \frac{1}{{\int_{X^{K-1}}^{X^{K}} \frac{dx}{A(x)}}} \sum_{j=1}^N h_{j} \frac{g(x^{K}_{j-\frac12})}{A(x^{K}_{j-\frac12})} \Big ) \\
&\qquad \qquad + \Big(-g(x^{K}_N) +\frac{1}{{\int_{X^{K}}^{X^{K+1}} \frac{dx}{A(x)}}} \sum_{j=1}^{N^{K+1}} h_{j+1} \frac{g(x^{K+1}_{j-\frac12})}{A(x^{K+1}_{j-\frac12})}      \Big) \\
&= \frac{1}{\int_{X^{K}}^{X^{K+1}} \frac{dx}{A(x)}} \sum_{j=1}^{N^{K+1}} h_{j+1} \frac{g(x^{K+1}_{j-\frac12})}{A(x^{K+1}_{j-\frac12})} - \frac{1}{\int_{X^{K-1}}^{X^{K}} \frac{dx}{A(x)}} \sum_{j=1}^{N^K}  h_{j} \frac{g(x^{K}_{j-\frac12})}{A(x^{K}_{j-\frac12})} \\
&= \frac{1}{\sum_{j=1}^{N^{K+1}} \Big( \big[ A^h \big]_{M_{K}+j,M_{K}+j-1}\Big)^{-1}} \sum_{j=1}^{N^{K+1}} \frac{g(x_{M_{K}+j-\frac12})}{\big[ A^h \big]_{M_{K}+j,M_{K}+j-1}} \\
& \qquad \qquad -\frac{1}{\sum_{j=1}^{N^K} \Big( \big[ A^h \big]_{M_{K-1}+j,M_{K-1}+j-1}\Big)^{-1} }  \sum_{j=1}^{N^K} \frac{g(x_{M_{K-1}+j-\frac12})}{\big[ A^h \big]_{M_{K-1}+j,M_{K-1}+j-1}}.
\end{aligned}
\end{equation*}
For $K=0$ and $K=N_H$,
\begin{equation*}
\begin{aligned}
&b^{H}_{0} = -g(x_{0})+\frac{1}{\sum_{j=1}^{N^1} \Big( \big[ A^h \big]_{j,j-1}\Big)^{-1} }  \sum_{j=1}^{N^1} \frac{g(x_{j-\frac12})}{\big[ A^h \big]_{j,j-1}}, \\
&b^{H}_{N_{H}} = g(x_{N_h})-\frac{1}{\sum_{j=1}^{N^{N^{H}}} \Big( \big[ A^h \big]_{M_{N^{H}-1}+j,M_{N^{H}-1}+j-1}\Big)^{-1} }  \sum_{j=1}^{N^{N^{H}}} \frac{g(x_{M_{N^{H}-1}+j-\frac12})}{\big[ A^h \big]_{M_{N^{H}-1}+j,M_{N^{H}-1}+j-1}}.
\end{aligned}
\end{equation*}
Since we know $g(x_{0}),g(x_{\frac12}),g(x_{\frac32}), \cdots,
g(x_{N_{h}-\frac12}), g(x_{N_{h}})$ from the micro--scale right hand
side vectors as given in \eqref{eq:rhs_micro4},
the macro--scale right hand side vector can be computed
from the micro--scale matrix system.

\subsection{Multiscale solution} Now we get the macro--scale matrix system using algebraic structure of micro--scale matrix system. We solve this macro--scale system to get the multiscale solution
\begin{equation*}
u_{ms}^{H} = \sum_{K=0}^{N^{H}} \eta^{K} \Psi^{K}.
\end{equation*}
Note that for $1 \leq K \leq N^{H}-1,$
\begin{equation*}
\Psi^{K}(x) = 
\left\{
\begin{aligned}
& \Psi^{K}_{+}(x) \quad \quad \text{ for all } x \in I^{K}, \\
& \Psi^{K+1}_{-}(x) \quad \text{ for all } x \in I^{K+1}.
\end{aligned}
\right.
\end{equation*}
For $K=0$ and $K=N_H$,
\begin{equation*}
\Psi^{0}(x) = \Psi^{0}_{-}(x) \quad \text{ and } \quad \Psi^{N_H}(x) = \Psi^{N_H}_{+}(x). 
\end{equation*}
We already know the exact form of $\Psi^{K}_{\pm}$ from \eqref{eq:formularPsiKhpm}. Using this form, we can compute the value of $\Psi^{K}$ at each micro node with the micro--scale stiffness matrix elements. Note that for $ 1 \leq l \leq N^{K}$,
\begin{eqnarray*}
\int_{X^{K-1}}^{x^{K}_{l}}\frac{1}{A(s)} \; ds &=& \sum_{j=1}^{l} \int_{x^{K}_{j-1}}^{x^{K}_{j}} \frac{dx}{A(x)}
\approx \sum_{j=1}^{l} \frac{h_{j}}{A(x^{K}_{j-\frac12})}
\\ &=& \sum_{j=1}^{l} \frac{h_{j}}{A_{M_{K-1}+j-\frac12}}
=-\sum_{j=1}^{l} \Big(\big[A^h \big]_{M_{K-1}+j,M_{K-1}+j-1}  \Big)^{-1}.
\end{eqnarray*}
Thus
\begin{equation*}
\Psi^{K}_{+}(x^{K}_{l}) = \frac{\sum_{j=1}^{l} \Big(\big[A^h \big]_{M_{K-1}+j,M_{K-1}+j-1} \Big)^{-1}}{\sum_{j=1}^{N^K} \Big( \big[ A^h \big]_{M_{K-1}+j,M_{K-1}+j-1}\Big)^{-1} }.
\end{equation*}
Similarly,
\begin{equation*}
\Psi^{K}_{-}(x^{K}_{l}) = \frac{\sum_{j=l+1}^{N^{K}} \Big(\big[A^h \big]_{M_{K-1}+j,M_{K-1}+j-1} \Big)^{-1}}{\sum_{j=1}^{N^K} \Big( \big[ A^h \big]_{M_{K-1}+j,M_{K-1}+j-1}\Big)^{-1} }.
\end{equation*}
Now we can easily compute the value of $u_{ms}^{H}$ at every micro node.

%
%
\section{Numerical examples}
In this section, we investigate some numerical examples to see the
convergence behavior of our scheme. We take $N^{h} = 2^{10}$.
The micro--scale solution $u^{h}$ is used as a reference solution to
compute error. Since $u_{ms}^{H}$ and $u^{h}$ are piecewise linear
functions, we have, for all $1 \leq j \leq N_{h}$,
\begin{equation*}
\begin{aligned}
\frac{d u_{ms}^{H}}{dx}\Big|_{[x_{j-1},x_{j}]} &= \frac{u_{ms}^{H}(x_{j})-u_{ms}^{H}(x_{j-1})}{h_{j}} = \frac{\Delta u_{ms}^{H}(j)}{h_j}, \\
\frac{d u^{h}}{dx}\Big|_{[x_{j-1},x_{j}]} &= \frac{u^{h}(x_{j})-u^{h}(x_{j-1})}{h_{j}} = \frac{\Delta u^{h}(j)}{h_j}.
\end{aligned}
\end{equation*}
The energy--norm error $E(u_{ms}^{H})$ of $u_{ms}^{H}$ is computed by
\begin{equation*}
\begin{aligned}
E(u_{ms}^{H}) &= \int_{0}^{1} A(x) \big[ \frac{d}{dx} \big( u_{ms}^{H}(x) \big) \big]^{2} \; dx \Big\}^{\frac12}
=\Big\{ \sum_{j=1}^{N_{h}} \int_{x_{j-1}}^{x^{j}} A(x) \Big( \frac{\Delta u_{ms}^{H}(j)}{h_{j}} \Big)^{2}
\Big\}^{\frac12} \\
&\approx \Big\{ \sum_{j=1}^{N_{h}} \frac{1}{h_{j}} A(x_{j-\frac12}) \big(\Delta u_{ms}^{H}(j) \big)^{2} \Big\}^{\frac12}
= \Big\{ \sum_{j=1}^{N_{h}} \big[ A^h \big]_{j,j-1} \big(\Delta
u_{ms}^{H}(j) \big)^{2} \Big\}^{\frac12}.
\end{aligned}
\end{equation*} 
The energy--norm $E(u_{h})$ of $u_{h}$ is obtained similarly. In the
following examples, we compute the relative energy--norm error:
\begin{equation*}
e_{energy}^{H} = \frac{E(u_{ms}^{H} - u^{h})}{E(u^{h})}.
\end{equation*}
\subsection{Known Coefficient Case}
\begin{example}\label{ex1}
Consider the equation \eqref{eq:1d-ell} with $\O=(0,1)$, $f=-1$ and
$A(x) = \frac{2}{3}(1+x)(1+\cos(\frac{2\pi x}{\epsilon})^2),$
$\epsilon=\frac1{10}.$ The homogenized solution is given by
\begin{equation*}
u_{Hom}(x) = \frac{3}{2\sqrt{2}}\Big(x-\frac{\log{(1+x)}}{\log{2}}\Big).
\end{equation*}
First we consider uniform micro--scale meshes. In Table~\ref{tbl:ex1_unif},
$e_{energy}^H$ and $e_{L^{2}}^{H}$ denote the relative energy-norm and
$L^2$ errors, respectively, of our AMS solutions. The errors are optimal.
\begin{table}
\begin{center}
\begin{tabular}{|c||cc|cc|cc|}
\hline
$N^{H}$&$e_{energy}^H$ & Order &$e_{L^{2}}^{H}$ & Order \\ \hline
2      & 5.00E-01      &       & 3.24E-01       &       \\ \hline               
4      & 2.50E-01      &  1.00 & 6.74E-01       &  2.26 \\ \hline
8	   & 1.26E-01      &  0.99 & 1.63E-02       &  2.05 \\ \hline
16     & 6.27E-02      &  1.00 & 4.07E-03       &  2.00 \\ \hline
32     & 3.08E-02      &  1.02 & 9.73E-04       &  2.06 \\ \hline
64     & 1.55E-02      &  0.99 & 2.51E-04       &  1.95 \\ \hline
\end{tabular}
\end{center}
\caption{Energy-norm and
$L^2$ errors and reduction rates of \exref{ex1} with uniform micro--scale meshes.}
\label{tbl:ex1_unif}
\end{table}
In Fig. \ref{fig:ex1-basis} first multiscale basis functions
$\Psi^{1}$ of \exref{ex1} are shown with $N^H=2,4,8,16,$ where uniform
micro--scale meshes are used. In Fig. \ref{fig:ex1_sol}, the red circle line and the blue
dashed-cross line denote the micro--scale solution $u^{h}$ and
multiscale solution $u_{ms}^{H}$, respectively with $N^H=2,4,8,16,64$.
We observe that the
multiscale solution $u_{ms}^H$ converges to the micro--scale solution
$u^{h}$ as the size $N^{H}$ becomes larger.

\begin{figure}
\begin{center}
\includegraphics[width=0.7\columnwidth]{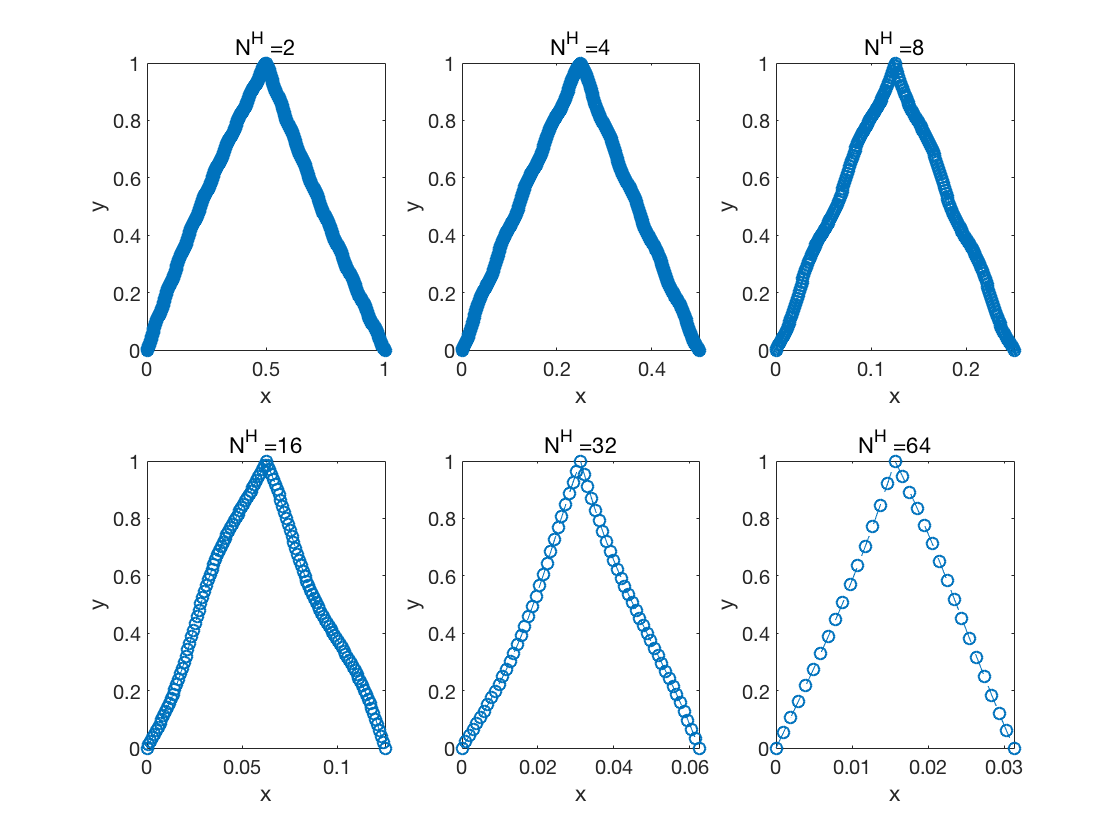}
\caption{First macro--scale basis function $\Psi^{1}$ of \exref{ex1}
  with uniform micro--scale mesh.}
\label{fig:ex1-basis}
\end{center}
\end{figure}

\begin{figure}
\begin{center}
  \includegraphics[width=0.45\columnwidth]{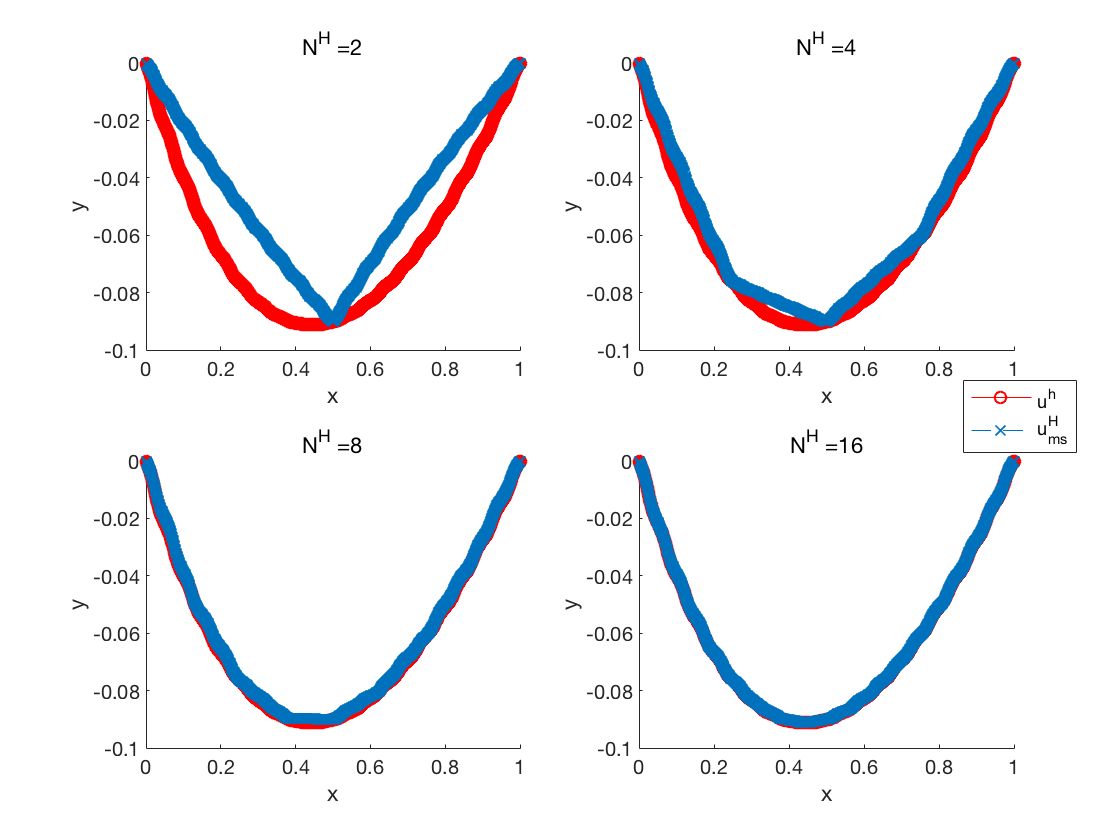}
  \includegraphics[width=0.45\columnwidth]{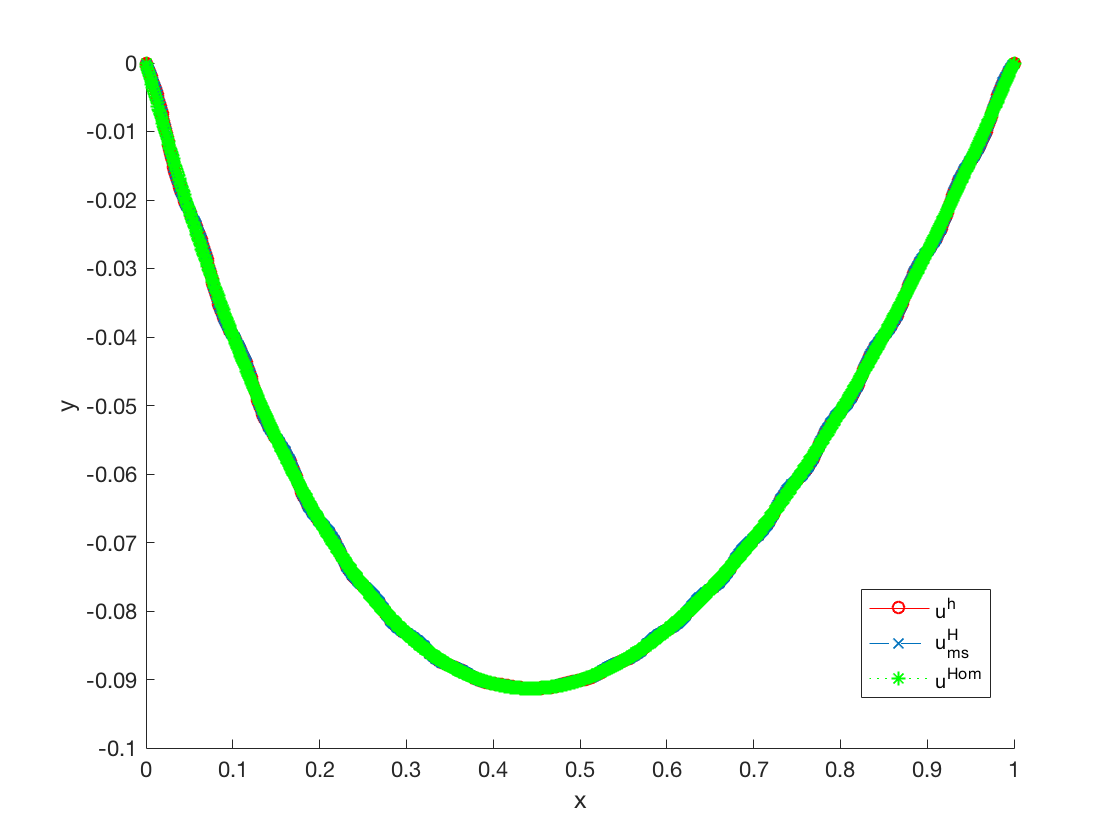}
\caption{Graph of solutions of \exref{ex1} with uniform micro--scale
  meshes: (L) with $N^H=2,4,8,16;$ (R) with $N^{H}=64.$ o
The red-circled and the blue dashed-cross lines denote the
micro--scale solution $u^{h}$ and multiscale solution $u_{ms}^{H}$,
respectively. The green dotted star line denotes the homogenized solution.
}
\label{fig:ex1_sol}
\end{center}
\end{figure}

Next, consider non-uniform micro--scale meshes. Let $y_{0}=0$ and
\begin{equation*}
y_{j+1} = y_{j} + 2\times \frac{\rand}{N_{h}} \text{  for  }  1 \leq j \leq N_{h},
\end{equation*}
where the MATLAB \it{rand} function is used. Take $x = \frac{y}{y_{N_{h}}}$ as
a micro--scale mesh so that $x_{N_h}=1.$

\begin{table}
\begin{center}
\begin{tabular}{|c||cc|cc|cc|}
\hline
$N^{H}$&$e_{energy}^H$ & Order &$e_{L^{2}}^{H}$ & Order \\ \hline
2      & 5.01E-01      &       & 3.22E-01       &       \\ \hline               
4      & 2.50E-01      &  1.00 & 6.63E-02       &  2.28 \\ \hline
8	   & 1.26E-01      &  0.99 & 1.60E-02       &  2.05 \\ \hline
16     & 6.28E-02      &  1.00 & 4.05E-03       &  1.99 \\ \hline
32     & 3.12E-02      &  1.01 & 9.95E-04       &  2.02 \\ \hline
64     & 1.60E-02      &  0.96 & 2.67E-04       &  1.90 \\ \hline
\end{tabular}
\end{center}
\caption{Energy-norm and
$L^2$ errors  and reduction rates of \exref{ex1} with non-uniform micro--scale meshes.}
\label{tbl:ex1_nonunif}
\end{table}

\begin{figure}
\begin{center}
\includegraphics[width=0.7\columnwidth]{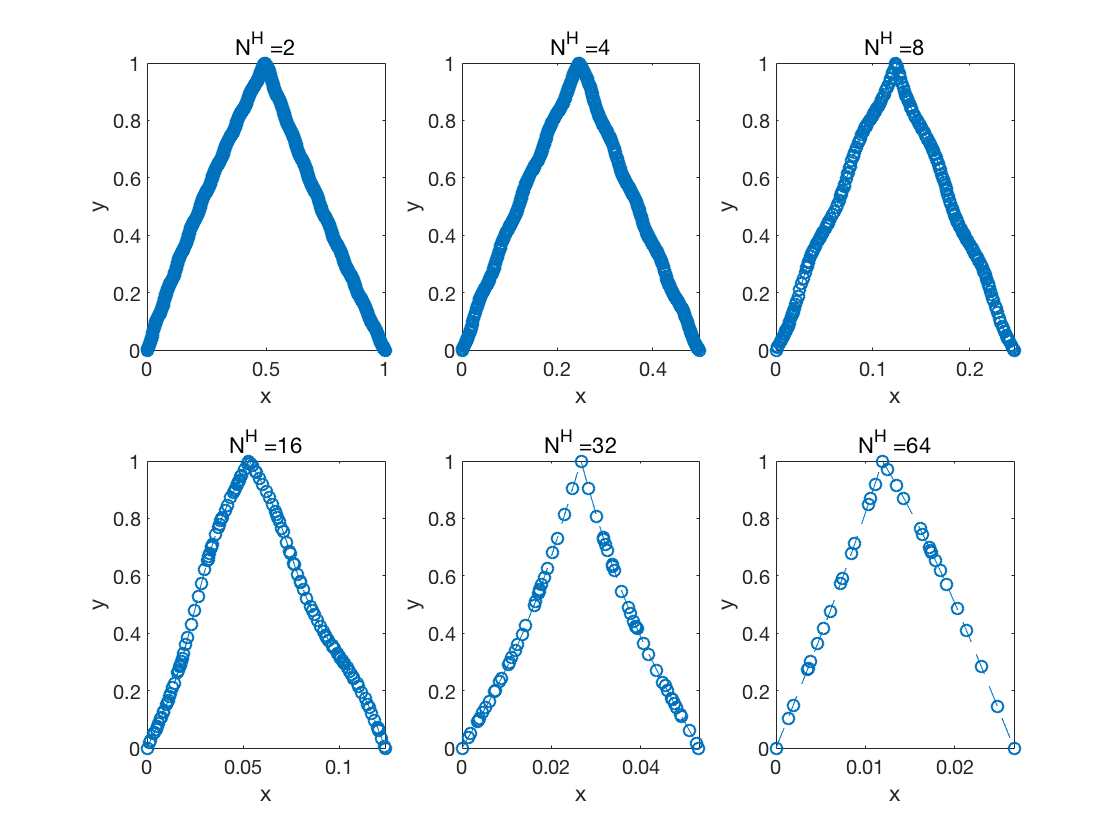}
\caption{First macro--scale basis function $\Psi^{1}$ of \exref{ex1} with non-uniform micro--scale mesh.}
\end{center}
\end{figure}\label{fig:ex1-nonuni-basis}

\begin{figure}
\begin{center}
  \includegraphics[width=0.45\columnwidth]{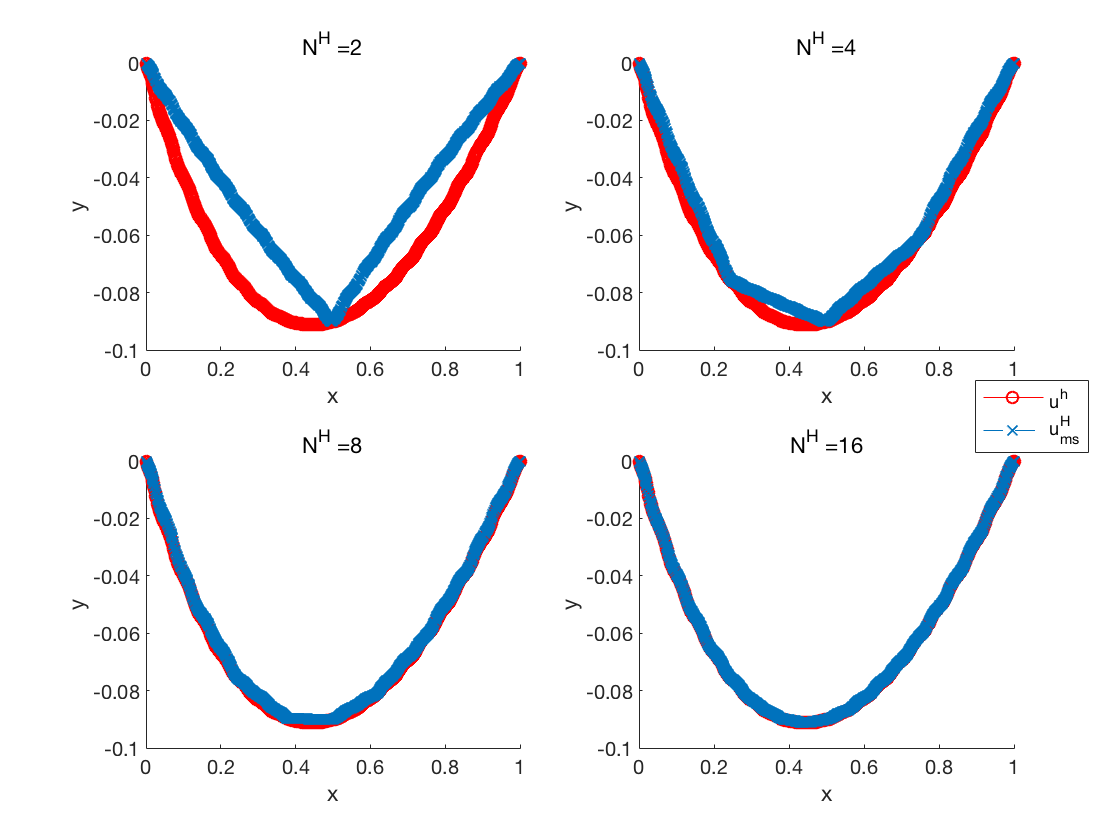}
  \includegraphics[width=0.45\columnwidth]{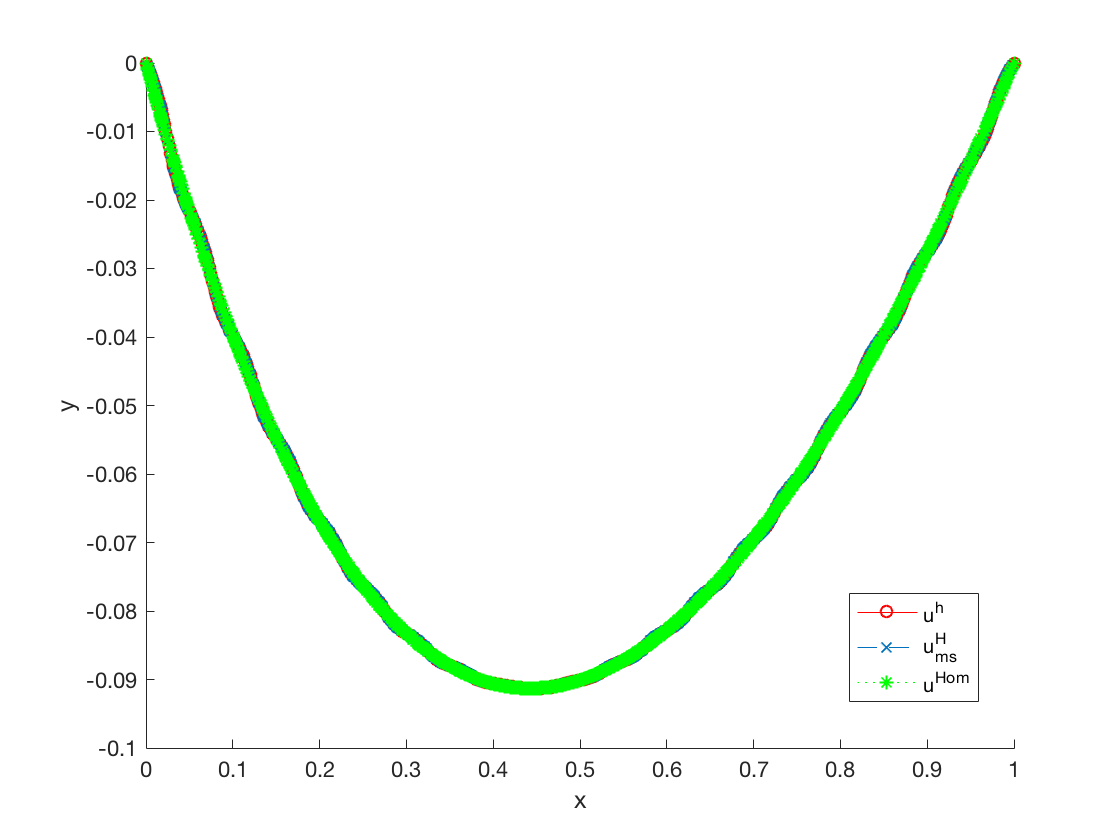}
\caption{
Graph of solutions of \exref{ex1} with non--uniform micro--scale
  meshes: (L) with $N^H=2,4,8,16;$ (R) with $N^{H}=64.$
The red-circled and the blue dashed-cross lines denote the
micro--scale solution $u^{h}$ and multiscale solution $u_{ms}^{H}$,
respectively. The green dotted star line denotes the homogenized solution.
}
\end{center}
\end{figure}
Tab. \ref{tbl:ex1_nonunif} shows the energy and $L^2$--errors for
nonuniform meshes and Figs. 4.3 and 4.4
show the graph of
solutions and first multiscale basis functions.
We observe that that the multiscale solution $u_{ms}^H$ converges to
the micro--scale solution $u^{h}$ as the size $N^{H}$ becomes larger.

\end{example}
\subsection{Random Coefficient Case}
We consider \eqref{eq:1d-ell} with $\O=(0,1)$, and the values of $A$ and $f$
are given randomly. That is, we consider the
situation that only the micro--scale matrix system is known without
any information on the exact form of $A$ and $f$. Also we can further
assume that we do not know the geometric information of micro--scale mesh. In our simulation, the micro--scale right hand side is given by
\begin{equation*}
b^{h}_{j} = \rand \text{ for } 0 \leq j \leq N_{h}.
\end{equation*}
For the micro--scale stiffness matrix we need to infer
off--diagonal elements. To construct a symmetric positive
definite tridiagonal matrix $A^h$ we use \eqref{eq:micro_stiff_diag}
and \eqref{eq:micro_stiff_off-diag}. We remark that the off--diagonal
element of $A^h$ corresponds to the average of coefficient matrix in
each micro interval.

We consider four Examples 4.2--4.5.
The cases of periodic coefficients with fixed amplitude and growing
amplitude are considered in Examples 4.2 and 4.3, respectively, while
those of non--periodic case in Examples 4.4 and 4.5, respectively.
Error tables are given in Tabs. 4.3--4.6, and the graphs of
coefficient, macro--scale basis functions, and solution graphs are shown
Figs. 4.5--4.12, respectively. Notice that in particular for Example
4.4, the random coefficient leads to discontinuous macro basis
functions $\Psi^1$ (see \eqref{eq:macro-basis-problem1}, \eqref{eq:macro-basis-problem2}) in Fig. 4.9 and the numerical approximation by the AMS
method in Fig. 4.10. Similar results are shown in Figs. 4.11--4.12. In
all the cases the approach by the AMS recovers the macro--scale
solutions based on the information on the micro--scale linear systems only.
These solutions recover very well the macro--scale solutions obtained by the
standar GMsFEM which are available when the explicit information on
the coefficient $A(x)$ in micro scale and the source function $f(x)$ are known.
\begin{example}\label{ex2}[Periodic behavior keeping its initial amplitude]
In this example the off--diagonal elements of $A^h$ are defined by
\begin{equation*}
\big[ A^h \big]_{j,j-1} = \big( 2 - 0.5\cos(\rand \times 17.7 j \pi) \big) \text{ for } 1\leq j \leq N_{h}.
\end{equation*}

\begin{figure}
\begin{center}
  \includegraphics[width=0.45\columnwidth]{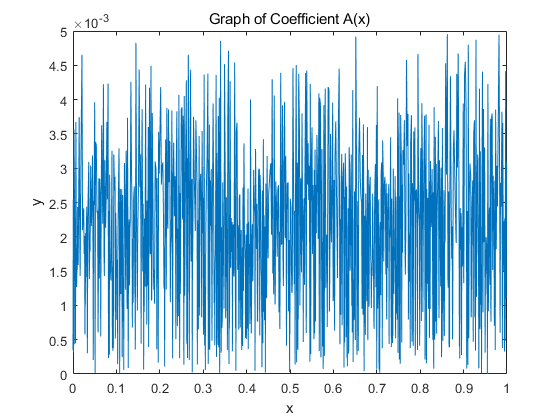}
  \includegraphics[width=0.45\columnwidth]{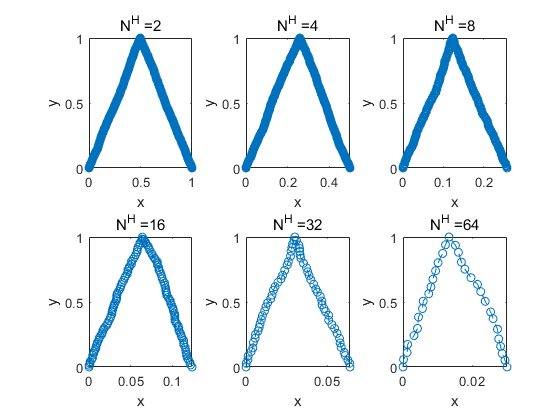}
\caption{
Graph of coefficient $A(x)$ of \exref{ex2} and first multiscale basis functions $\Psi^{1}$.}
\end{center}
\end{figure}

\begin{table}
\begin{center}
\begin{tabular}{|c||cc|cc|cc|}
\hline
$N^{H}$&$e_{energy}^H$ & Order &$e_{L^{2}}^{H}$ & Order \\ \hline
2      & 5.04E-01      &       & 3.20E-01       &       \\ \hline               
4      & 2.51E-01      &  1.01 & 6.62E-02       &  2.27 \\ \hline
8	   & 1.27E-01      &  0.99 & 1.60E-02       &  2.05 \\ \hline
16     & 6.30E-02      &  1.01 & 3.96E-03       &  2.02 \\ \hline
32     & 3.17E-02      &  0.99 & 9.95E-04       &  1.99 \\ \hline
64     & 1.59E-02      &  1.00 & 2.48E-04       &  2.01 \\ \hline
\end{tabular}
\end{center}
\caption{Energy-norm and $L^2$ errors  and reduction rates of \exref{ex2}.}
\label{tbl:ex2}
\end{table}

\begin{figure}
\begin{center}
  \includegraphics[width=0.45\columnwidth]{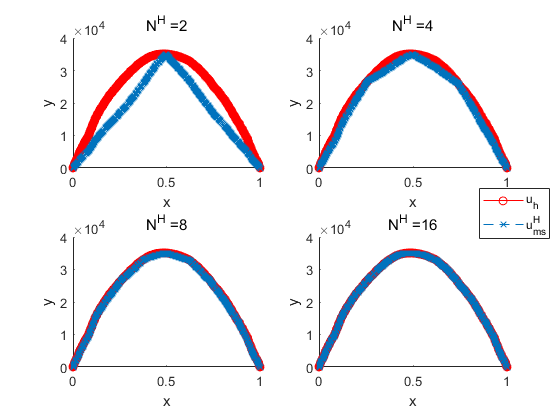}
  \includegraphics[width=0.45\columnwidth]{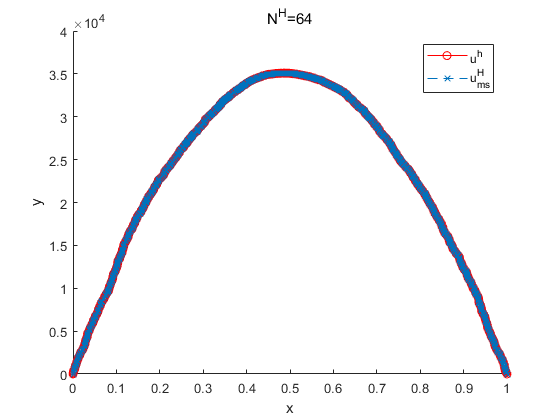}
\caption{Graph of solutions of \exref{ex2} with non--uniform micro--scale
  meshes: (L) with $N^H=2,4,8,16;$ (R) with $N^{H}=64.$ 
The red-circled and the blue dashed-cross lines denote the
micro--scale solution $u^{h}$ and multiscale solution $u_{ms}^{H}$,
respectively. The green dotted star line denotes the homogenized solution.
}
\end{center}
\end{figure}
\end{example}

\begin{example}\label{ex3}[Periodic behavior with growing amplitude]
Define the off--diagonal element of  $A^h$ by
\begin{equation}
\big[ A^h \big]_{j,j-1} = j\big( 2 - 0.5\cos(\rand \times 17.7 j \pi) \big) \text{ for } 1\leq j \leq N_{h}.
\end{equation}

\begin{figure}
\begin{center}
  \includegraphics[width=0.45\columnwidth]{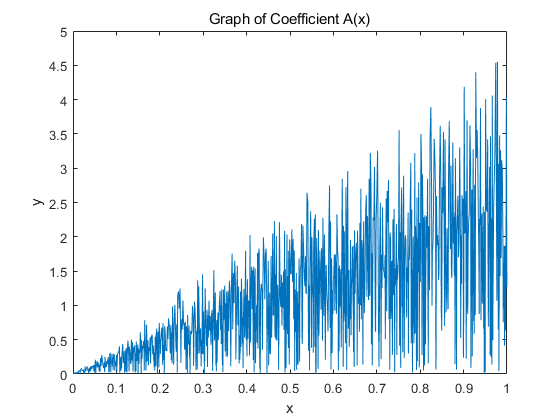}
  \includegraphics[width=0.45\columnwidth]{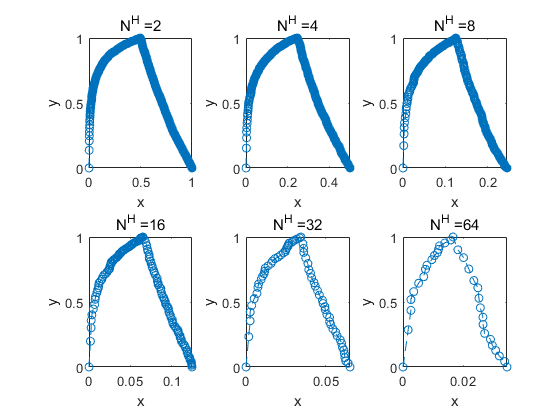}
\caption{
Graph of coefficient $A(x)$ of \exref{ex3} and first multiscale basis functions $\Psi^{1}$.
}
\end{center}
\end{figure}

\begin{table}
\begin{center}
\begin{tabular}{|c||cc|cc|cc|}
\hline
$N^{H}$&$e_{energy}^H$ & Order &$e_{L^{2}}^{H}$ & Order \\ \hline
2      & 5.27E-01      &       & 4.83E-01       &       \\ \hline               
4      & 2.79E-01      &  0.92 & 1.27E-01       &  1.92 \\ \hline
8	   & 1.44E-01      &  0.96 & 3.95E-02       &  1.69 \\ \hline
16     & 7.54E-02      &  0.93 & 1.29E-02       &  1.61 \\ \hline
32     & 3.93E-02      &  0.94 & 4.37E-03       &  1.56 \\ \hline
64     & 2.19E-02      &  0.85 & 1.60E-03       &  1.44 \\ \hline
\end{tabular}
\end{center}
\caption{Energy-norm and $L^2$ errors  and reduction rates of \exref{ex3}.}
\label{tbl:ex3}
\end{table}

\begin{figure}
\begin{center}
  \includegraphics[width=0.45\columnwidth]{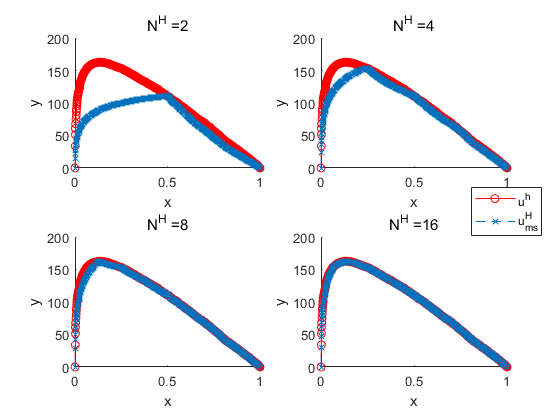}
  \includegraphics[width=0.45\columnwidth]{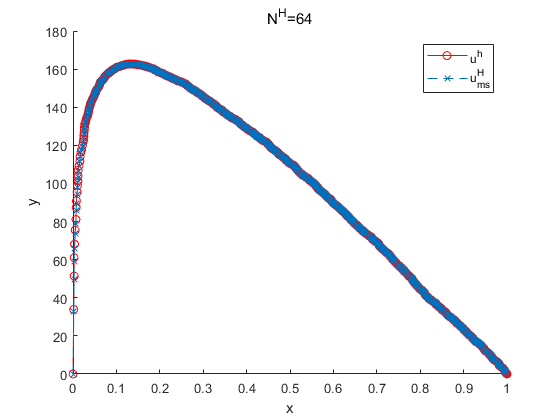}
  \caption{
Graph of solutions of \exref{ex3} with non--uniform micro--scale
  meshes: (L) with $N^H=2,4,8,16;$ (R) with $N^{H}=64.$ 
The red-circled and the blue dashed-cross lines denote the micro--scale solution $u^{h}$ and
multiscale solution $u_{ms}^{H}$, respectively. 
  }
\end{center}
\end{figure}

\end{example}

\begin{example}\label{ex4}[Non-periodic behavior keeping its initial amplitude]
Define the off--diagonal element of  $A^h$ by
\begin{equation}
\big[ A^h \big]_{j,j-1} = \rand \text{ for } 1\leq j \leq N_{h}.
\end{equation}

\begin{figure}
\begin{center}
  \includegraphics[width=0.45\columnwidth]{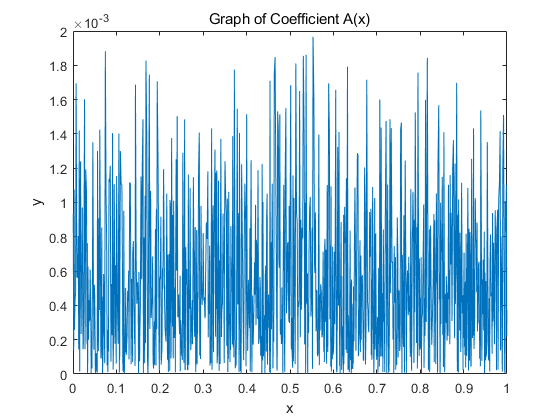}
  \includegraphics[width=0.45\columnwidth]{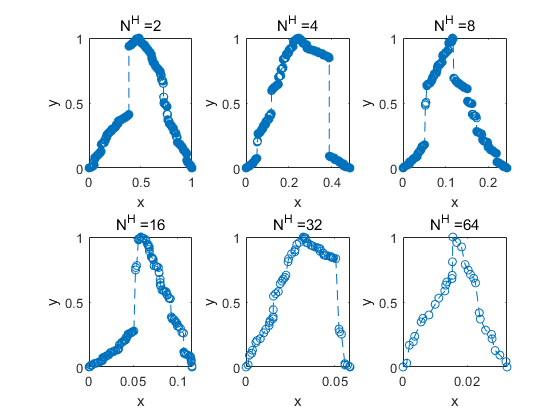}
\caption{
Graph of coefficient $A(x)$ of \exref{ex4} and first multiscale basis functions $\Psi^{1}$.}
\end{center}
\end{figure}

\begin{table}
\begin{center}
\begin{tabular}{|c||cc|cc|cc|}
\hline
$N^{H}$&$e_{energy}^H$ & Order &$e_{L^{2}}^{H}$ & Order \\ \hline
2      & 5.16E-01      &       & 3.68E-01       &       \\ \hline               
4      & 2.37E-01      &  1.12 & 5.96E-02       &  2.63 \\ \hline
8	   & 1.13E-01      &  1.07 & 1.37E-02       &  2.12 \\ \hline
16     & 5.43E-02      &  1.06 & 3.57E-03       &  1.94 \\ \hline
32     & 2.48E-02      &  1.13 & 6.82E-04       &  2.39 \\ \hline
64     & 1.30E-02      &  0.94 & 2.29E-04       &  1.58 \\ \hline
\end{tabular}
\end{center}
\caption{Energy-norm and $L^2$ errors  and reduction rates of \exref{ex4}.}
\label{tbl:ex4}
\end{table}

\begin{figure}
\begin{center}
  \includegraphics[width=0.45\columnwidth]{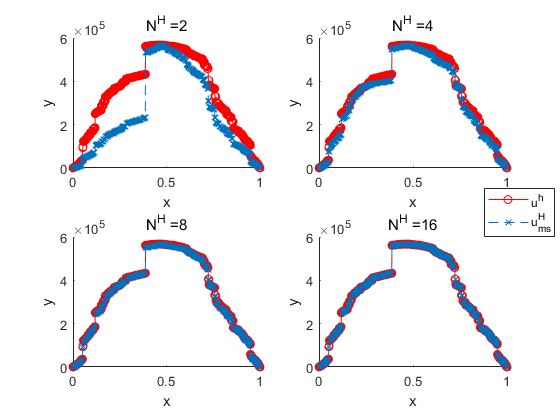}
  \includegraphics[width=0.45\columnwidth]{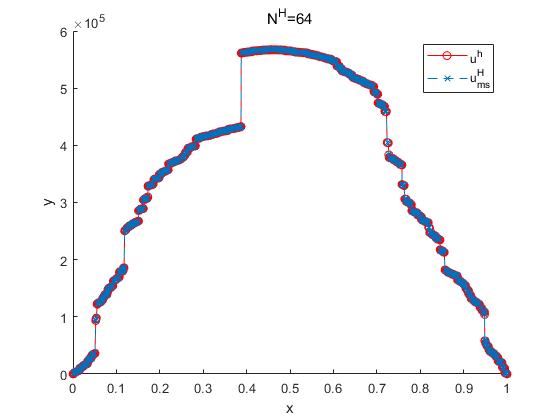}
  \caption{
Graph of solutions of \exref{ex4} with non--uniform micro--scale
  meshes: (L) with $N^H=2,4,8,16;$ (R) with $N^{H}=64.$
The red-circled and the blue
dashed-cross lines denote the micro--scale solution $u^{h}$ and
multiscale solution $u_{ms}^{H}$, respectively. 
}
\end{center}
\end{figure}

\end{example}

\begin{example}\label{ex5}[Non--periodic behavior with growing amplitude]
Define the off--diagonal element of  $A^h$ by
\begin{equation}
\big[ A^h \big]_{j,j-1} = j*\rand \text{ for } 1\leq j \leq N_{h}.
\end{equation}

\begin{figure}
\begin{center}
  \includegraphics[width=0.45\columnwidth]{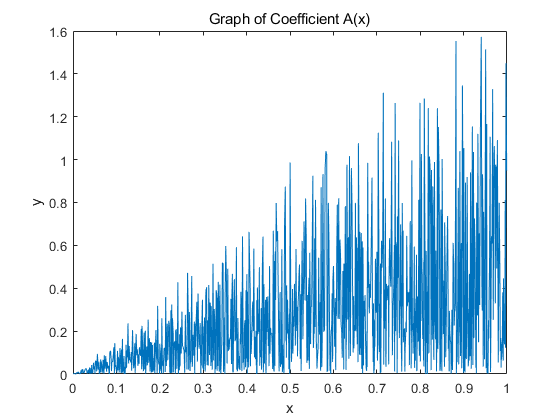}
  \includegraphics[width=0.45\columnwidth]{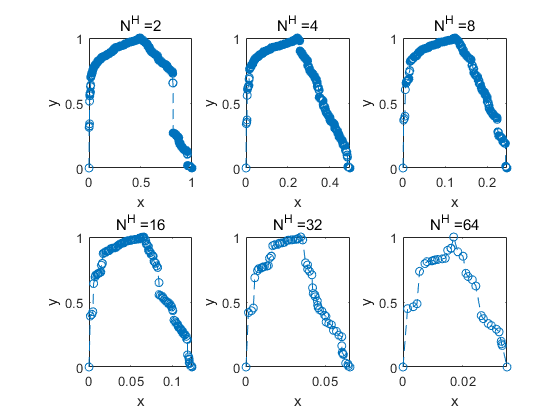}
\caption{Graph of coefficient $A(x)$ of \exref{ex5} and first multiscale basis functions $\Psi^{1}$.}
\end{center}
\end{figure}

\begin{table}
\begin{center}
\begin{tabular}{|c||cc|cc|cc|}
\hline
$N^{H}$&$e_{energy}^H$ & Order &$e_{L^{2}}^{H}$ & Order \\ \hline
2      & 4.00E-01      &       & 2.39E-01       &       \\ \hline               
4      & 2.03E-01      &  0.98 & 6.34E-02       &  1.91 \\ \hline
8	   & 9.29E-01      &  1.13 & 1.57E-02       &  2.02 \\ \hline
16     & 4.74E-02      &  0.97 & 5.16E-03       &  1.60 \\ \hline
32     & 2.55E-02      &  0.90 & 2.29E-04       &  1.17 \\ \hline
64     & 1.56E-02      &  0.71 & 1.37E-04       &  0.74 \\ \hline
\end{tabular}
\end{center}
\caption{Energy-norm and $L^2$ errors  and reduction rates of \exref{ex5}.}
\label{tbl:ex5}
\end{table}

\begin{figure}
\begin{center}
  \includegraphics[width=0.45\columnwidth]{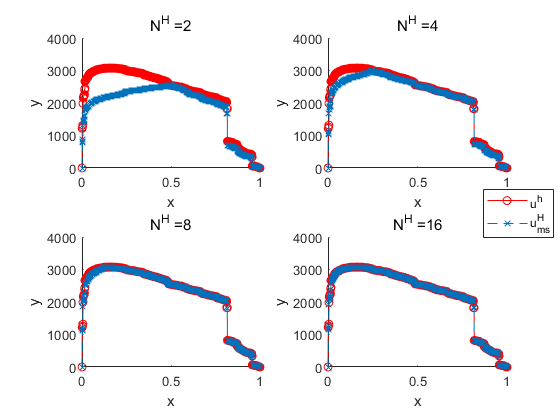}
  \includegraphics[width=0.45\columnwidth]{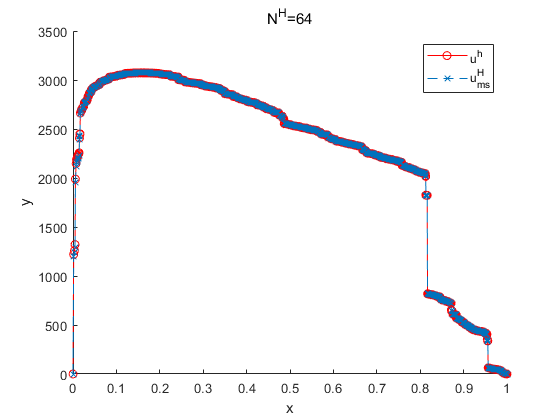}
  \caption{
  Graph of solutions of \exref{ex5} with non--uniform micro--scale
  meshes: (L) with $N^H=2,4,8,16;$ (R) with $N^{H}=64.$
The red-circled and the blue
dashed-cross lines denote the micro--scale solution $u^{h}$ and
multiscale solution $u_{ms}^{H}$, respectively. 
}
\end{center}
\end{figure}
\end{example}

%
%
\section{Conclusion}
In this paper we propose to a method to compute macro--scale solutions
based on the information on the linear system generated by a micro
scale elliptic equation by the finite element method. We do not assume
that
the function form for the coefficient and source term in micro scale
are not known. But we infer them from the coefficients of the given
linear system. We explain the detailed procedure how the 
We tested several numerical examples which confirm the approach
proposed in this paper works well for various cases.
\section*{ Acknowledgment} 
DS was supported in part by National Research Foundation of
Korea (NRF-2017R1A2B3012506 and NRF-2015M3C4A7065662). 

\normalsize
\bibliographystyle{abbrv}
\bibliography{ams}

\end{document}